\documentclass[11pt,a4paper]{amsart}

\usepackage[T1]{fontenc}
\usepackage[cyr]{aeguill}
\usepackage[english]{babel}
\usepackage[latin1]{inputenc}
\usepackage{fullpage}
\usepackage[dvips]{graphics}
\usepackage{color}
\usepackage{url}
\usepackage{stmaryrd}
\usepackage{verbatim}

\usepackage{amsmath}
\usepackage{amsthm}
\usepackage{amssymb}
\usepackage{bbm}
\usepackage{delarray}
\usepackage{enumerate}
\usepackage{longtable}
\usepackage[all,poly,necula]{xy}
\usepackage{ifthen}

\DeclareSymbolFont{bbdold}{U}{bbold}{m}{n}
\DeclareSymbolFontAlphabet{\mathbbd}{bbdold}

\theoremstyle{plain}
\newtheorem{thm}{Theorem}

\newtheorem{cor}[thm]{Corollary}

\newcounter{thmnct}

\theoremstyle{definition}

\theoremstyle{remark}

\makeatletter
\@addtoreset{footnote}{page}
\def\@fnsymbol#1{\ensuremath{\ifcase#1\or \dagger\or \ddagger\or {\dagger\dagger}\or {\ddagger\ddagger} \else{*}\fi}}

\makeatother

\DeclareMathOperator{\SL}{SL}

\DeclareMathOperator{\id}{Id}

\DeclareMathOperator{\Hom}{Hom}

\DeclareMathOperator{\md}{mod}

\DeclareMathOperator{\irr}{irr}

\renewcommand{\geq}{\geqslant}

\newcommand{\Z}{\mathbb{Z}}

\newcommand{\C}{\mathbb{C}}

\renewcommand{\epsilon}{\varepsilon}

\newcommand{\inter}{\cap}

\newcommand{\restr}[1]{\!|_{#1}}

\newcommand{\tens}{\otimes}

\renewcommand{\phi}{\varphi}

\renewcommand{\bar}[1]{\overline{#1}}

\newcommand{\forloop}[5][1] { \setcounter{#2}{#3} \ifthenelse{#4} { #5 \addtocounter{#2}{#1} \forloop[#1]{#2}{\value{#2}}{#4}{#5} }{ } } 
\newcounter{i}
\newcounter{j}
\newcommand{\dets}[4]{
    \displaystyle \left|\;
    \begin{matrix}
      \forloop{i}{1}{\value{i} < 7}{
       \forloop{j}{1}{\value{j} < 7}{
        \ifthenelse{\value{i} < #1 \or \value{i} > #2 \or \value{j} < #3 \or \value{j} > #4}{\ifthenelse{\value{i}<\value{j}}{\cdot}{\ifthenelse{\value{i} = \value{j}}{1}{}}}{\bullet} \ifthenelse{\value{j}<6}{&}{}
       } \ifthenelse{\value{i}<6}{\\}{}
      }
    \end{matrix}
    \;\right|^{\vphantom{M^{M^M}}}_{\vphantom{M_{M_M}}}}

\renewcommand{\tilde}[1]{\widetilde{#1}}

\makeatletter
\newcommand{\gl}{\begin@guill}
\newcommand{\gr}{\end@guill}
\makeatother

\let\oldmarginpar\marginpar
\renewcommand\marginpar[1]{\-\oldmarginpar[\raggedleft\footnotesize \color{blue}{#1}]%
{\raggedright\footnotesize \color{blue}{#1}}}

\begin{document}
\setlength{\parindent}{5mm}
\title{Skew group algebras of path algebras\\and preprojective algebras}
\author{Laurent Demonet}
\address{LMNO, Universit\'e de Caen, Esplanade de la Paix, 14000 Caen}
\email{Laurent.Demonet@normalesup.org}

\date{}

\begin{abstract}
 We compute explicitly up to Morita-equivalence the skew group algebra of a finite group acting on the path algebra of a quiver and the skew group algebra of a finite group acting on a preprojective algebra. These results generalize previous results of Reiten and Riedtmann \cite{ReRi85} for a cyclic group acting on the path algebra of a quiver and of Reiten and Van den Bergh \cite[proposition 2.13]{ReVa89} for a finite subgroup of $\SL(\C X \oplus \C Y)$ acting on $\C[X, Y]$.
\end{abstract}

\maketitle

\section{Introduction and main results} \label{intro}

Let $k$ be an algebraically closed field and $G$ be a finite group such that the characteristic of $k$ does not divide the cardinality of $G$. If $\Lambda$ is a $k$-algebra and if $G$ acts on the right on $\Lambda$, the action being denoted exponentially, the skew group algebra of $\Lambda$ under the action of $G$ is by definition the $k$-algebra whose underlying $k$-vector space is $\Lambda \tens_k k[G]$ and whose multiplication is linearly generated by $(a \tens g) (a' \tens g') = aa'^{g^{-1}} \tens gg'$ for all $a, a' \in \Lambda$ and $g, g' \in G$ (see \cite{ReRi85}). It will be denoted by $\Lambda G$. Identifying $k[G]$ and $\Lambda$ with subalgebras of $\Lambda G$, an alternative definition is
$$\Lambda G = \langle \Lambda, k[G] \,|\, \forall (g, a) \in G \times \Lambda, g^{-1} a g = a^g \rangle_{k\text{-alg}}$$

Let now $Q = (I, A)$ be a quiver where $I$ denotes the set of vertices and $A$ the set of arrows. Consider an action of $G$ on the path algebra $kQ$ permuting the set of primitive idempotents $\{e_i \,|\, i \in I\}$ and stabilizing the vector space spanned by the arrows. Note that this is more general than an action coming from an action of $G$ on $Q$ since an arrow may be sent to a linear combination of arrows. We now define a new quiver $Q_G$. We first need some notation.

Let $\tilde I$ be a set of representatives of the classes of $I$ under the action of $G$. For $i \in I$, let $G_i$ denote the subgroup of $G$ stabilizing $e_i$, let $i_\circ \in \tilde I$ be the representative of the class of $i$ and let $\kappa_i \in G$ be such that $i_\circ^{\kappa_i} = i$.

For $(i, j) \in \tilde I^2$, $G$ acts on $O_i \times O_j$ where $O_i$ and $O_j$ are the orbits of $i$ and $j$ under the action of $G$. A set of representatives of the classes of this action will be denoted by $F_{ij}$. 

For $i, j \in I$, define $M_{ij} \subset k Q$ to be the vector space spanned by the arrows from $i$ to $j$. We regard $M_{ij}$ as a right $k[G_i \inter G_j]$-module by restricting the action of $G$.

The quiver $Q_G$ has vertex set 
 $$I_G = \bigcup_{i \in \tilde I} \{i\} \times \irr(G_i)$$
 where $\irr(G_i)$ is a set of representatives of isomorphism classes of irreducible representations of $G_i$. The set of arrows of $Q_G$ from $(i, \rho)$ to $(j, \sigma)$ is a basis of
 $$\bigoplus_{(i', j') \in F_{ij}} \Hom_{\md k[G_{i'} \inter G_{j'}]}\left((\rho \cdot \kappa_{i'}) \restr{G_{i'} \inter G_{j'}}, (\sigma \cdot \kappa_{j'}) \restr{G_{i'} \inter G_{j'}} \tens_k M_{i'j'} \right) $$
 where the representation $\rho \cdot \kappa_{i'}$ of $G_{i'}$ is the same as $\rho$ as a vector space, and $(\rho \cdot \kappa_{i'})_{g} = \rho_{\kappa_{i'} g \kappa_{i'}^{-1}}$ for $g \in G_{i'} = \kappa_{i'}^{-1} G_i \kappa_{i'}$. Table \ref{t1} gives two examples of quivers $Q_G$. A detailed example is also computed in section \ref{ex}.
  
 \begin{table}[b] \caption{Examples of computations of $Q_G$} \label{t1}
  \newcommand{\lmp}[1] {\begin{minipage}{6cm}#1\vspace*{0cm}\end{minipage}}
  \begin{tabular}{|c|c|c|}
    \hline
    $Q$ & $G$ & $Q_G$ \\
    \hline
    \lmp{
    $$\xymatrix@R=.49cm@C=.49cm{
     1 \ar[r] & 2 \ar[r] & \dots \ar[r] & n-1 \ar[rd] &  \\
      & & & & n \\
     1' \ar[r] & 2' \ar[r] & \dots \ar[r] & (n-1)' \ar[ru] & 
    }$$}
    &
    $\Z/2\Z$
    &
    \lmp{$$\xymatrix@R=.49cm@C=.49cm{
     & & & & n_+ \\
     1 \ar[r] & 2 \ar[r] & \dots \ar[r] & n-1 \ar[ur] \ar[dr] &  \\
     & & & & n_- \\
    }$$} \\ 
    \hline
    \lmp{$$\xymatrix@R=1.49cm@C=1.49cm{
     1 \ar@(dl,dr)[]_\beta \ar@(ul,ur)[]^\alpha
    }$$}
    &
    \begin{minipage}{2.5cm}
    \begin{center}
    $G \subset \SL\left(\C \alpha \oplus \C \beta \right)$ \\    
    type $A_n$
    \end{center}
    \end{minipage}
    &
    \lmp{$$\xymatrix@R=.49cm@C=.49cm{
       & 0 \ar@/^/[dr] \ar@/^/[dl] & \\
       1 \ar@/^/[d] \ar@/^/[ur] & & n \ar@/^/[d] \ar@/^/[ul] \\
       2 \ar@/^/[r] \ar@/^/[u] & \dots \ar@/^/[r] \ar@/^/[l] & n-1 \ar@/^/[u] \ar@/^/[l] 
    }$$} \\ 
    \hline
  \end{tabular}
 \end{table}

We can now state the two main results of this paper and two corollaries :

\begin{thm} \label{th1}
 There is an equivalence of categories
 $$\md k\left(Q_G\right) \simeq \md \left(kQ\right)G.$$
\end{thm}
 Theorem \ref{th1} was proved by Reiten and Riedtmann in \cite[\S 2]{ReRi85} for cyclic groups.

The following theorem deals with the case of preprojective algebras. The definition of the preprojective algebra $\Lambda_Q$ of a quiver $Q$ is recalled in section \ref{prth2}.

\begin{thm} \label{th2}
  If $G$ acts on $k \bar Q$, where $\bar Q$ is the double quiver of $Q$, by permuting the primitive idempotents $e_i$ and stabilizing the linear subspace of $k \bar Q$ spanned by the arrows, and if, for all $g \in G$, $r^g = r$ where $r$ is the preprojective relation of this quiver, then $\left(\bar Q\right)_G$ is of the form $\bar Q'$ for some quiver $Q'$ and $(\Lambda_Q)G$ is Morita equivalent to $\Lambda_{Q'}$.
\end{thm}

One can always extend an action on $kQ$ to an action on $k \bar Q$ and this yields :
 
\begin{cor} \label{c1}
 The action of $G$ on a path algebra $k Q$ permuting the primitive idempotents and stabilizing the linear subspace of $k \bar Q$ spanned by the arrows induces naturally an action of $G$ on $k \bar Q$ and $\left(\bar Q\right)_G$ is isomorphic to the double quiver of $Q_G$. Moreover, there is an equivalence of categories
  $$\md \Lambda_{Q_G} \simeq \md \Lambda_Q G.$$
\end{cor}

Theorem \ref{th2} and corollary \ref{c1} will be used in \cite{De} for constructing $2$-Calabi-Yau categorifications of skew-symmetrizable cluster algebras. Another corollary, which is an easy consequence of the definition of the McKay graph, is :

\begin{cor} \label{c2}
 Let $Q$ be the quiver 
 $$
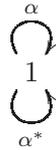
\xymatrix@R=1.49cm@C=1.49cm{
     1 \ar@(dl,dr)[]_{\alpha^*} \ar@(ul,ur)[]^\alpha
    }$$
 and $G$ a finite subgroup of $\SL(\C \alpha \oplus \C \alpha^*)$. Then 
 \begin{enumerate}
  \item There is an identification of $Q_G$ with a double quiver $\bar Q'$ such that the non oriented underlying graph of $Q'$ is isomorphic to the affine Dynkin diagram corresponding to $G$ through the McKay correspondence.
  \item We have $kQ / (\alpha \alpha^* - \alpha^* \alpha) \simeq k[\alpha, \alpha^*]$ and there is an equivalence of categories
  $$\md \left(k[\alpha, \alpha^*] G\right) \simeq \md \Lambda_{Q'}.$$
 \end{enumerate}
\end{cor}

 Corollary \ref{c2} was proved by a geometrical method in \cite[proof of proposition 2.13]{ReVa89} (see also \cite[theorem 0.1]{CrHo98}).

\section{An example} \label{ex}

Suppose that $k = \C$ and that $Q$ is the following quiver 
$$
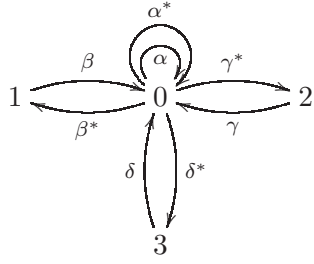
\xymatrix@R=1.49cm@C=1.49cm{
     1 \ar@/^/[r]^\beta & \ar@/^/[l]^{\beta^*} 0 \ar@`{p+(-7,4),p+(0,15),p+(7,4)}[]^{\alpha^*} \ar@(ul,ur)[]_\alpha \ar@/^/[r]^{\gamma^*} \ar@/^/[d]^{\delta^*} & 2 \ar@/^/[l]^{\gamma} \\
      & 3 \ar@/^/[u]^{\delta} &
    }$$
Let also
$$G = \langle a,b \,|\, a^3 = b^2, b^4 = 1, aba = b\rangle$$
be the binary dihedral group of order $12$. One lets $G$ act on $kQ$ by :
 \begin{center} \vspace{.1cm}
 \begin{tabular}{|c|c|c|c|c|c|c|c|c|c|c|c|c|}
  \cline{2-13}
   \multicolumn{1}{c|}{}& $e_0$ & $e_1$ & $e_2$ & $e_3$ & $\alpha$ & $\alpha^*$ & $\beta$ & $\beta^*$ & $\gamma$ & $\gamma^*$ & $\delta$ & $\delta^*$ \\ 
  \hline
   $a$ & $e_0$ & $e_2$ & $e_3$ & $e_1$ & $\zeta^{-1} \alpha$ & $\zeta \alpha^*$ & $\gamma$ & $\gamma^*$ & $\delta$ & $\delta^*$ & $\beta$ & $\beta^*$ \\ 
   \hline
   $b$ & $e_0$ & $e_1$ & $e_3$ & $e_2$ & $\alpha^*$ & $-\alpha$ & $-\beta$ & $-\beta^*$ & $-\delta$ & $-\delta^*$ & $-\gamma$ & $-\gamma^*$ \\
   \hline
 \end{tabular} \vspace{.1cm}
 \end{center}
 where $\zeta$ is a primitive sixth root of unity.

 Using the notation of the introduction, one can choose $\tilde I = \{0, 1\}$, $\kappa_0 = \kappa_1 = 1$, $\kappa_2 = a$, $\kappa_3=a^2$. One has $G_0 = G$, $G_1 = \langle b \rangle \simeq \Z/4\Z$, $G_2 = \langle ba \rangle \simeq \Z/4\Z$, $G_3 = \langle ab \rangle \simeq \Z/4\Z$. One can also choose $F_{0,0} = \{(0,0)\}$, $F_{0,1} = \{(0,1)\}$, $F_{1,0} = \{(1,0)\}$ and $F_{1,1} = \{(1,1), (1,2), (2,1)\}$. 
 
 The irreducible representations of $\Z/4\Z$ will be denoted by $\theta_\alpha$ where $\alpha \in \{i, -1, -i, 1\}$ is the scalar action of a specified generator ($b$, $ba$ or $ab$ when $\Z/4\Z$ is realized as $G_1$, $G_2$ or $G_3$). The group $G$ has six irreducible representations : four of degree $1$ of the form $a \mapsto \alpha^2$, $b \mapsto \alpha$ for each $\alpha \in \{i, -1, -i, 1\}$, which will be denoted by $\lambda_\alpha$, and two of degree $2$ :
 $$\rho : a \mapsto \left( \begin{matrix} \zeta^{-1} & 0 \\ 0 & \zeta \end{matrix} \right) \quad b \mapsto \left( \begin{matrix} 0 & -1 \\ 1 & 0 \end{matrix} \right)$$
 and
 $$\sigma : a \mapsto \left( \begin{matrix} \zeta^{-2} & 0 \\ 0 & \zeta^2 \end{matrix} \right) \quad b \mapsto \left( \begin{matrix} 0 & 1 \\ 1 & 0 \end{matrix} \right)$$
 One checks easily that $\lambda_i \tens \rho \simeq \sigma$, $\rho \tens \rho \simeq \sigma \oplus \lambda_1 \oplus \lambda_{-1}$, $\rho \tens \sigma \simeq \rho \oplus \lambda_i \oplus \lambda_{-i}$, $\sigma \tens \sigma \simeq \sigma \oplus \lambda_1 \oplus \lambda_{-1}$. The other product formulas are deduced from these. One computes $M_{0,0} = \rho$, $M_{0,1} = M_{(1,0)} = \lambda_{-1}$ and $M_{(1,1)} = M_{(1,2)} = M_{(2,1)} = 0$.
 The vertices of $Q_G$ are then $0_{i}$, $0_{-1}$, $0_{-i}$, $0_{1}$, $0_{\rho}$, $0_{\sigma}$, $1_{-1}$ and $1_{1}$ where we write $0_\alpha = (0, \lambda_\alpha)$ and $1_\alpha = (1, \theta_\alpha)$ for simplicity. One has
 $$\Hom_G(\rho, \sigma \tens \rho) \simeq \Hom_G(\rho, \rho \oplus \lambda_i \oplus \lambda_{-i}) \simeq \C$$
 and therefore there is one arrow from $0_\rho$ to $0_\sigma$,
 $$\Hom_G(\lambda_1, \sigma \tens \rho) \simeq 0$$
 and therefore there is no arrow from $0_1$ to $0_\sigma$,
 $$\Hom_G(\lambda_i, \sigma \tens \rho) \simeq \C$$
 and therefore there is one arrow from $0_i$ to $0_\sigma$,
 $$\Hom_{\Z/4\Z}(\theta_i, \sigma\restr{\Z/4\Z} \tens \theta_{-1}) \simeq \Hom_{\Z/4\Z}(\theta_i, \theta_{-1} \oplus \theta_1) \simeq 0$$
 and therefore there is no arrow from $1_i$ to $0_\sigma$,
 $$\Hom_{\Z/4\Z}(\theta_1, \sigma\restr{\Z/4\Z} \tens \theta_{-1}) \simeq \Hom_{\Z/4\Z}(\theta_1, \theta_{-1} \oplus \theta_1) \simeq \C$$
 and therefore there is one arrow from $1_1$ to $0_\sigma$,
 $$\Hom_{\Z/4\Z}(\theta_1, \lambda_{-1} \restr{\Z/4\Z} \tens \theta_{-1}) \simeq \C$$
 and therefore there is one arrow from $1_1$ to $0_{-1}$. All the other computations can be done in the same way. Finally, $Q_G$ is the following quiver :
 $$
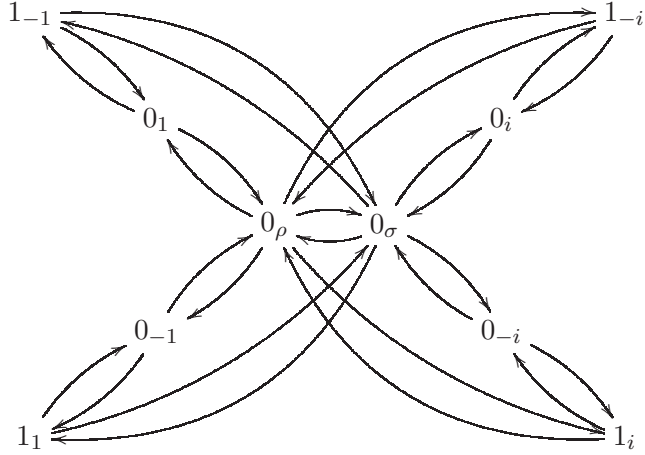
\xymatrix{ 
  1_{-1} \ar@/^/[dr] \ar@/^1cm/[ddrrr] & & & & & 1_{-i} \ar@/^/[dl] \ar@/_.5cm/[ddlll] \\
  & 0_1 \ar@/^/[ul] \ar@/^/[dr] & & & 0_i \ar@/^/[ur] \ar@/^/[dl] \\
  & & 0_\rho \ar@/^/[ul] \ar@/^/[dl] \ar@/^/[r] \ar@/^1cm/[uurrr] \ar@/_.5cm/[ddrrr] & 0_\sigma \ar@/^/[l] \ar@/^/[ur] \ar@/^/[dr] \ar@/_.5cm/[uulll] \ar@/^1cm/[ddlll] & & \\
  & 0_{-1} \ar@/^/[ur] \ar@/^/[dl] & & & 0_{-i} \ar@/^/[ul] \ar@/^/[dr] \\
  1_1 \ar@/^/[ur] \ar@/_.5cm/[uurrr]& & & & & 1_i \ar@/^/[ul] \ar@/^1cm/[uulll]
 }$$
 where one can remark that the full subgraph having vertices $\{0_i, 0_{-1}, 0_{-i}, 0_1, 0_\rho, 0_\sigma\}$ is the affine Dynkin diagram corresponding to $G$ in the McKay correspondence, as expected. Hence $\md (\C Q)G \simeq \md \C Q_G$. Moreover, it is easy to check that the preprojective relation $\alpha \alpha^* - \alpha^* \alpha + \beta \beta^* - \beta^* \beta + \gamma \gamma^* - \gamma^* \gamma + \delta \delta^* - \delta^* \delta$ is stable under the action of $G$ and therefore there is an equivalence of Morita between $\Lambda_{Q_1} G$ and $\Lambda_{Q_2}$ where $\bar Q_1 = Q$, $\bar Q_2 = Q_G$, and $\Lambda_{Q_1}$, $\Lambda_{Q_2}$ are the preprojective algebras of $Q_1$ and $Q_2$.

\section{Proofs of the main propositions} \label{prv}

\subsection{Proof of theorem \ref{th1}} \label{prth1}
We retain the notation of section \ref{intro}. Thus $Q$ is a quiver, $G$ acts on $kQ$ by stabilizing the set of primitive idempotents corresponding to vertices and the vector space spanned by the arrows, $\tilde I$ is a fixed set of representatives of the $G$-orbits of $I$, etc.. Let $R$ be the subalgebra of $kQ$ generated by the primitive idempotents and $M \subset kQ$ be the linear subspace spanned by the arrows, seen as an $R$-bimodule. Define $T_0 = R$ and for every positive integer $n$, $T_n = T_{n-1} \tens_R M$. Then, recall that the tensor algebra $T(R, M)$ is $\bigoplus_{i \geq 0} T_i$ endowed with the canonical product. It is clear that $k Q$ is canonically isomorphic to $T(R, M)$ on which the action of $G$ is graded.  

As $G$ stabilizes $R$ and $M$, one can define the skew-group algebra $RG$ which is a subalgebra of $(kQ) G$. Thus, $MG = M \tens_k k[G]$ is a sub-$RG$-bimodule of $(kQ)G$ and one gets easily a canonical isomorphism $(kQ) G \simeq T(R,M) G \simeq T(RG, MG)$ which maps $(m_1 \tens m_2 \tens \dots \tens m_n) \tens g$ to $\left(m_1 \tens 1\right) \tens \left(m_2 \tens 1\right) \tens \dots \tens \left(m_{n-1} \tens 1\right) \tens \left(m_n \tens g\right) = \left(m_1 \tens 1\right) \tens \left(m_2 \tens 1\right) \tens \dots \tens \left(m_{n-1} \tens g\right) \tens \left(m_n^g \tens 1\right) = \dots = \left(m_1 \tens 1\right) \tens \left(m_2 \tens g\right) \tens \dots \tens \left(m_{n-1}^g \tens 1\right) \tens \left(m_n^g \tens 1\right) = \left(m_1 \tens g\right) \tens \left(m_2^g \tens 1\right) \tens \dots \tens \left(m_{n-1}^g \tens 1\right) \tens \left(m_n^g \tens 1\right)$.

Let now $e = \sum_{i \in \tilde I} e_i \in R \subset RG$ which is an idempotent. Then, according to \cite[proposition 1.6]{ReRi85} and its proof,
\begin{align*}
 \md RG & \simeq \md e(RG)e \\
 N & \mapsto eN \\
 (RG)e \tens_{e(RG)e} N' & \mapsfrom N'
\end{align*}
is a Morita equivalence between $RG$ and $e(RG)e$. Hence, it is a classical and easy fact that there is an equivalence of categories
\begin{align*}
 \md (kQ)G \simeq \md T(RG, MG) & \simeq \md T(e(RG)e, e(MG)e) \\
 N & \mapsto eN \\
 (RG)e \tens_{e(RG)e} N' & \mapsfrom N'
\end{align*}
using the isomorphism $MG \simeq (RG)e \tens_{e(RG)e} e(MG)e \tens_{e(RG)e} e(RG)$. One has $e(RG)e \simeq \prod_{i \in \tilde I} k[G_i]$ where, for each $i \in \tilde I$, $G_i$ is the stabilizer of $e_i$. As $G_i$ is semi-simple, one can fix, for each $i \in \tilde I$ and $\rho \in \irr(G_i)$, $\tilde e_{i\rho}$ to be a primitive idempotent of $k[G_i]$ corresponding to $\rho$. Let $\tilde e = \sum_{i \in \tilde I} \sum_{\rho \in \irr(G_i)} \tilde e_{i \rho}$ which satisfies $e \tilde e e = \tilde e$. Then we have a Morita equivalence between $e(RG)e$ and $\tilde e e (RG) e \tilde e = \tilde e (RG) \tilde e$, and, as before, between $T(e(RG)e, e(MG)e)$ and $T(\tilde e(RG)\tilde e, \tilde e (MG)\tilde e)$. Moreover, using the notation $I_G$ of the introduction, $\tilde e(RG)\tilde e \simeq \prod_{(i, \rho) \in I_G} k \tilde e_{i \rho}$ and therefore, it is enough to compute $\tilde e_{j\sigma} (MG) \tilde e_{i \rho} = \tilde e_{j \sigma} e_j (MG) e_i \tilde e_{i \rho}$ for each $(i, \rho)$ and $(j, \sigma)$ in $I_G$ to finish the proof of theorem \ref{th1}. Remark now that 
 $$e_j (MG) e_i = \sum_{(i',j') \in O_i \times O_j} G_j \kappa_{j'} M_{i'j'} \kappa_{i'}^{-1} G_i = \bigoplus_{(i',j') \in F_{ij}} G_j \kappa_{j'} M_{i'j'} \kappa_{i'}^{-1} G_i$$
 and therefore, if one denotes $G_{i'j'} = G_{i'} \inter G_{j'}$ for every $i', j' \in I$,
 \begin{align*}\tilde e_{j\sigma} e_j (MG) e_i \tilde e_{i \rho} \simeq \bigoplus_{(i',j') \in F_{ij}} &\Hom_k \left(k, \tilde e_{j\sigma} G_j \kappa_{j'} M_{i'j'} \kappa_{i'}^{-1} G_i \tilde e_{i \rho}\right) \\
   \simeq \bigoplus_{(i',j') \in F_{ij}} &\Hom_{k[G_i]} \left(\rho, \sigma \tens_{k[G_j]} G_j \kappa_{j'} M_{i'j'} \kappa_{i'}^{-1} G_i \right) \\
   \simeq \bigoplus_{(i',j') \in F_{ij}} &\Hom_{k[G_{i'}]} \left(\rho \cdot \kappa_{i'}, (\sigma \cdot \kappa_{j'}) \tens_{k[G_{j'}]} G_{j'} M_{i'j'} G_{i'} \right) \\
   \simeq \bigoplus_{(i',j') \in F_{ij}} &\Hom_{k[G_{i'}]} \left(\rho \cdot \kappa_{i'}, (\sigma \cdot \kappa_{j'}) \tens_{k[G_{j'}]} k[G_{j'}] \right. \\ &\left. \tens_{k[G_{i'j'}]} G_{i'j'} M_{i'j'}  G_{i'j'} \tens_{k[G_{i'j'}]} k[G_{i'}]\right)
 \end{align*}
 and, because of the relations defining $(kQ) G$, the multiplication in $(kQ) G$ induces an isomorphism of $(G_{i'j'}, G_{i'j'})$-bimodules $k[G_{i'j'}]\tens_k M_{i'j'} \simeq G_{i'j'} M_{i'j'} G_{i'j'}$. Note that the action of $(G_{i'j'}, G_{i'j'})$ on $k[G_{i'j'}] \tens_k M_{i'j'}$ is defined here by $g (v \tens m) h = gvh \tens h^{-1} m h = gvh \tens m^h$ (recall that we use the multiplication notation for the product in $(kQ)G$ and the exponential notation for the action of $G$ on $kQ$ fixed at the beginning). Therefore 
 \begin{align*}\tilde e_{j\sigma} e_j (MG) e_i \tilde e_{i \rho} \simeq \bigoplus_{(i',j') \in F_{ij}} &\Hom_{k[G_{i'}]} \left(\rho \cdot \kappa_{i'}, (\sigma \cdot \kappa_{j'}) \vphantom{\tens_{k[G_{i'j'}]} \left(k[G_{i'j'}] \tens_k M_{i'j'}\right) \tens_{k[G_{i'j'}]}  k[G_{i'}]} \right. \\ & \left. \tens_{k[G_{i'j'}]} \left(k[G_{i'j'}] \tens_k M_{i'j'}\right) \tens_{k[G_{i'j'}]}  k[G_{i'}] \right)\\
   \simeq \bigoplus_{(i',j') \in F_{ij}} &\Hom_{k[G_{i'}]} \left(\rho \cdot \kappa_{i'}, \left((\sigma \cdot \kappa_{j'})\restr{G_{i'j'}} \tens_k  M_{i'j'}\right) \tens_{k[G_{i'j'}]} k[G_{i'}] \right)\\
   \simeq \bigoplus_{(i', j') \in F_{ij}} &\Hom_{k[G_{i'j'}]}\left((\rho \cdot \kappa_{i'}) \restr{G_{i'j'}}, (\sigma \cdot \kappa_{j'}) \restr{G_{i'j'}} \tens_k M_{i'j'} \right) 
 \end{align*}
 where the representation $\rho \cdot \kappa_{i'}$ of $G_{i'}$ is the same as $\rho$ as a vector space and, if $g \in G_{i'} = \kappa_{i'}^{-1} G_i \kappa_{i'}$, $(\rho \cdot \kappa_{i'})_g = \rho_{\kappa_{i'} g \kappa_{i'}^{-1}}$. It concludes the proof of theorem \ref{th1}. Note that we go from the penultimate line to the last one by the classical adjunction between induction and restriction of representations.

\subsection{Proof of theorem \ref{th2}} \label{prth2}

We retain the notation of section \ref{prth1}. Define the $R$-bimodule $\bar M = M \oplus M^*$ (the $R$-bimodule structure on $M^*$ is the natural one : $(afb)(m) = f(bma)$ for $a, b \in R$, $f \in M^*$ and $m \in M$). The tensor algebra $T(R, \bar M)$ is the path algebra of the double quiver $\bar Q$ of $Q$. The non-degenerate skew-symmetric bilinear form defined on $\bar M$ by $\langle m + f, m' + f'\rangle = f'(m) - f(m')$ where $m, m' \in M$ and $f, f' \in M^*$ satisfies, for every $a, b \in R$ and $m, n \in \bar M$, $\langle amb, n \rangle = \langle m, bna \rangle$. If $\{x_i\}_{i \in S}$ is a $k$-basis of $\bar M$, then denote by $\{x_i^*\}_{i \in S}$ its (left) dual basis for the bilinear form $\langle -, - \rangle$ (that is the one satisfying $\langle x_i^*, x_j \rangle = \delta_{ij}$ for every $i, j \in S$). Then the element $r = \sum_{i \in S} x_i \tens x_i^* \in \bar M \tens_R \bar M$ is independent of the choice of the basis. By definition, $r$ is the preprojective relation corresponding to $Q$ and the preprojective algebra $\Lambda_Q$ of $Q$ is $T(R, \bar M)/(r)$. For more details about preprojective algebras, see for example \cite{DlRi80} or \cite{Ri98}.

Suppose now that the group $G$ acts on $k \bar Q$ by stabilizing the set of primitive idempotents corresponding to vertices and the $k$-subspace spanned by the arrows. Then it stabilizes $r$ if and only if it stabilizes the bilinear form $\langle -, -\rangle$ (indeed, for $g \in G$, $r^g = r$ implies that $\{x_i^{*g}\}_{i \in S}$ is the left dual basis of $\{x_i^g\}_{i \in S}$ since for every $i \in S$, $(\id_{\bar M} \tens \langle -, - \rangle)(r^g \tens x_i^g) = (\id_{\bar M} \tens \langle -, - \rangle)(r \tens x_i^g ) = x_i^g$). Extend now $\langle -,-\rangle$ to a bilinear form on $\bar M G$ by setting
 $$\langle m \tens g, n \tens h \rangle = \left\{ \begin{array}{ll}
    \langle m, n^h \rangle & \text{if } gh = 1 \\
    0 & \text{else} \end{array} \right.$$
 which is clearly skew-symmetric and non-degenerate. Moreover, for $a, b \in RG$ and $m, n \in \bar M G$, one has easily $\langle amb, n \rangle = \langle m, bna \rangle$. If $\{x_i\}_{i \in S}$ is a basis of $\bar M$ and $\{x_i^*\}_{i \in S}$ is its left dual basis, then $\{x_i^{*g} \tens g^{-1}\}_{(i,g) \in S \times G}$ is the left dual basis of the basis $\{x_i \tens g\}_{(i,g) \in S \times G}$ of $\bar M G$. Hence, the preprojective relation $r_G$ corresponding to $\langle -, - \rangle$ in $\bar MG \tens_{RG} \bar MG$ is
  \begin{align*}r_G &= \sum_{(i, g) \in S \times G} \left(x_i \tens g\right) \tens \left(x_i^{*g} \tens g^{-1}\right) = \sum_{(i, g) \in S \times G} \left(x_i \tens 1\right)\left(1 \tens g\right) \tens \left(x_i^{*g} \tens g^{-1}\right) \\&= \sum_{(i, g) \in S \times G}\left(x_i \tens 1\right) \tens \left(1 \tens g\right)\left(x_i^{*g} \tens g^{-1}\right) = \sum_{(i, g) \in S \times G} \left(x_i \tens 1\right) \tens \left(x_i^*\tens 1\right) = \#G \times r\end{align*}
  where the preprojective relation $r$ of $T(R, \bar M)$ is mapped by the canonical inclusion from $\bar M \tens_R \bar M$ to $\bar MG \tens_{RG} \bar MG$. As $\# G$ is invertible in $k$, one gets 
  $$\left(T(R, \bar M) / (r)\right)G \simeq T(R G, \bar M G) / (r) = T(R G, \bar M G) / (r_G).$$
 Moreover, for $\langle -, - \rangle$, $\tilde e (\bar M G) \tilde e$ and $(1-\tilde e) (\bar M G) + (\bar M G) (1-\tilde e)$ are orthogonal supplementary subspaces. Thus, $\langle -, -\rangle$ restricts to a skew-symmetric non-degenerate bilinear form on $\tilde e (\bar M G) \tilde e$ which satisfies, for every $a,b \in \tilde e (RG) \tilde e$ and $m, n \in \tilde e (\bar M G) \tilde e$, $\langle amb, n \rangle = \langle m, bna \rangle$. By taking a basis of $\bar M G$ which is the union of a basis of $\tilde e (\bar M G) \tilde e$ and a basis of $(1-\tilde e) (\bar M G) + (\bar M G) (1-\tilde e)$, it is clear that the equivalence of categories of section \ref{prth1} restricts to an equivalence
 \begin{align*}
 \md \left((kQ)G / (r_G)\right) \simeq \md \left(T(RG, \bar M G)/(r_G)\right) & \simeq \md \left(T(\tilde e(RG) \tilde e, \tilde e(\bar M G) \tilde e)/(r_{\tilde e})\right) \\
 N & \mapsto \tilde eN \\
 (RG)\tilde e \tens_{\tilde e(RG)\tilde e} N' & \mapsfrom N'
\end{align*}
where $r_{\tilde e}$ is the preprojective relation in $T(\tilde e(RG) \tilde e, \tilde e(\bar M G) \tilde e) \simeq k \left(\bar Q\right)_G \simeq k \bar Q'$. It completes the proof of theorem \ref{th2} (it is enough to take a maximal isotropic subspace of $\tilde e(\bar M G) \tilde e$ to find arrows of $Q'$ which is of course non unique).

\subsection{Proof of corollary \ref{c1}} We retain the previous notation. If $G$ acts on $k Q$ by stabilizing the set of primitive idempotents and $M$, then its action can be extended to an action on $\bar M = M \oplus M^*$ using the contragredient representation on $M^*$. Then $\bar M G = M G \oplus M^* G$ as $RG$-bimodules. Moreover, $M G$ is clearly maximal isotropic for the bilinear form $\langle -,- \rangle$ extended on $\bar M G$ as before. Hence $\left(\bar Q\right)_G = \overline{Q_G}$ and it proves corollary \ref{c1}.

\section*{Acknowledgments}
The author would like to thank his PhD advisor Bernard Leclerc for his advices and corrections. He would also like to thank Christof Gei\ss, Bernhard Keller and Idun Reiten for interesting discussions and comments on the topic. He is also grateful to the referee for helpful suggestions for shortening and readability of this article.

\bibliographystyle{alphanum}
\bibliography{biblio}

\end{document}